\input psfig.sty
%

\newcount\m \newcount\n \newcount\p \newdimen\dim

\def\today{\ifcase\month\or
  January\or February\or March\or April\or May\or June\or
  July\or August\or September\or October\or November\or December\fi
  \space\number\day, \number\year}

\def\hoje{\number\day\ de
  \ifcase\month\or
  janeiro\or fevereiro\or mar{\c c}o\or abril\or maio\or junho\or
  julho\or agosto\or setembro\or outubro\or novembro\or dezembro\fi
  \ de \number\year}

\def\japtoday{\number\year.\twodigits\month.\twodigits\day}

\def\hours{\n=\time \divide\n 60
  \m=-\n \multiply\m 60 \advance\m \time
  \twodigits\n:\twodigits\m}

\def\twodigits#1{\ifnum #1<10 0\fi \number#1}

\def\datepages{\footline={{\hss\tenrm\folio\hss\hbox to 0pt{
\hidewidth\fiverm\japtoday}}}}

\tolerance=4000
\magnification=\magstep1
\parskip=4pt plus 6pt

\font\bigbf=cmbx12

\font\vbigbf=cmbx12 scaled \magstep1

\font\tencmr=cmr10

\font\sevencmr=cmr7
\font\fivecmr=cmr5
\newfam\cmrfam
\textfont\cmrfam=\tencmr
\scriptfont\cmrfam=\sevencmr
\scriptscriptfont\cmrfam=\fivecmr

\font\eusms=eusm7
\font\eusmf=eusm5
\newfam\eusmfam
\textfont\eusmfam=\eusms
\scriptfont\eusmfam=\eusmf
\scriptscriptfont\eusmfam=\eusmf

\font\tenbf=cmbx10
\font\eightbf=cmbx8
\font\sevenbf=cmbx7
\font\fivebf=cmbx5
\newfam\cmbxfam
\textfont\cmbxfam=\tenbf
\scriptfont\cmbxfam=\sevenbf
\scriptscriptfont\cmbxfam=\fivebf

\font\tenmsbm=msbm10
\font\sevenmsbm=msbm7
\font\fivemsbm=msbm5
\newfam\msbmfam
\textfont\msbmfam=\tenmsbm
\scriptfont\msbmfam=\sevenmsbm
\scriptscriptfont\msbmfam=\fivemsbm
\def\msbm{\fam\msbmfam}

\def\ZZ{{\msbm Z}}

\def\RR{{\msbm R}}
\def\CC{{\msbm C}}

\def\nobf{\noindent \bf}
\def\nobigbf{\noindent\bigbf}
\def\qed{{\hfill\vrule height3pt width6pt depth3pt\smallskip}}

\def\emptyset{/\kern -0.67em{\scriptstyle\bigcirc}}

\mathchardef\braceld="37A
\mathchardef\bracerd="37B

\def\m@th{\mathsurround=0pt }

\def\downbrakfill{$\m@th
\;\braceld\leaders\vrule\hfill\bracerd\;$}

\def\overbrak#1{\mathop{\vbox{\ialign{##\crcr\noalign{\kern3pt}
\downbrakfill\crcr\noalign{\kern3pt\nointerlineskip}
$\hfil\displaystyle{#1}\hfil$\crcr}}}\limits}

\def\sG{{\cal G}}

\def\sM{{\cal M}}
\def\sR{{\cal R}}
\def\taub{{\overline\tau}}
\def\tr{{\hbox{tr}}}
\def\sc{{\sqrt{5}}}



\centerline{\vbigbf Spectra of semi-regular polytopes}
\bigskip
\centerline{\it Nicolau C. Saldanha and Carlos Tomei
\footnote{}{\sevenrm Research supported by MCT and CNPq, Brazil.}}
\bigskip

\bigbreak

{

\narrower

{\noindent \bf Abstract:}
We compute the spectra of the adjacency matrices
of the semi-regular polytopes.
A few different techniques are employed:
the most sophisticated,
which relates the 1-skeleton of the polytope to a Cayley graph,
is based on methods akin to those of Lov\'asz and Babai ([L], [B]).
It turns out that the algebraic degree of the eigenvalues is
at most 5, achieved at two 3-dimensional solids.

}

\bigbreak
{\nobigbf Introduction}

\nobreak\bigskip\nobreak

Important information about a graph is conveyed by its spectrum,
i.e., the set of eigenvalues of its adjacency matrix;
an extensive bibliography can be found in [CDS].
Symmetries of a graph are helpful in computing the spectrum ([PS]).
If the group of symmetries acts transitively on vertices,
Lov\'asz ([L]) and Babai ([B]) show how to 
apply techniques from representation theory
of groups to reduce significantly the algebraic degree of the problem.

In [ST], we computed the spectra of all regular polytopes
(i.e., of the graph formed by their vertices and edges);
somewhat surprisingly, all eigenvalues have algebraic degree at most 3.
The complete list of regular polytopes is known since Schl\"afli
([Sc], [C]).
{\sl Semi-regular polytopes} are
polytopes with regular faces whose isometry group
acts transitively on vertices.
Trivial examples are regular polytopes and,
in three dimensions, prisms and antiprisms.
Nontrivial 3-dimensional semi-regular polytopes
are known as {\sl Archimedean solids}, despite the fact that
the writings by Archimedes on the subject are lost ([C]).
Kepler ([K]) wrote the first available list of Archimedean solids
with a proof its completeness.
The semi-regular polytopes in higher dimensions
were known to Gosset ([Go])
and were extensively studied by Coxeter ([C2], [C3]),
but only recently Blind and Blind ([BB]) showed
that Gosset's list is complete.

In this paper, we compute the spectra of all semi-regular polytopes.
Section 1 contains our results:
characteristic polynomials of the adjacency matrices
are completely factored in $\ZZ[\tau]$, $\tau = (1 + \sc)/2$.
It turns out that the algebraic degree of the eigenvalues is
at most 5, achieved at two 3-dimensional solids.
Section 2 is devoted to the study of some Archimedean solids
whose spectra can be obtained by taking into account
a few planes of symmetry.
Gosset polytopes are described in section 3:
since they have very large groups of symmetry,
the technique used in [ST] for regular polytopes is convenient here.
In section 4, we show how to apply group representations
to compute the spectrum of a Cayley graph
by Lov\'asz's, Babai's and our modified version of their methods.
Chung and Sternberg ([CS]) used similar ideas to compute
the spectra of some regular graphs and of the buckyball molecule,
which corresponds to a weighted version of
$\langle 5, 6, 6 \rangle$ (in the notation of the next section).
Group representation techniques
are used in section 5 to compute the spectra
of the remaining Archimedean solids.
In section 6 we consider three discrete subgroups of $S^3$,
the unit quaternionic sphere.
These come in handy in sections 7 and 8 where we compute the spectra
of the two missing semi-regular polytopes.
We tried to make this paper accessible to readers
with a scant knowledge of the semi-regular polytopes;
most difficulties should be resolved by consulting
the excellent book [C].

\bigbreak
{\nobigbf 1. Statement of results}

\nobreak\bigskip\nobreak

The semi-regular polytopes in dimension three are 
the five Platonic and the thirteen Archimedean solids,
besides prisms and antiprisms.
We shall denote the Archimedean solids
by the types of polygons surrounding each vertex.
Thus, $\langle 3,6,6\rangle$ is the solid
with two hexagons and a triangle at
each vertex, which is obtained by removing small
tetrahedra of edge 1 from each vertex
of a regular tetrahedron of edge 3.
In this notation, the Archimedean solids are
$\langle 3,6,6\rangle$, $\langle 3,8,8\rangle$,
$\langle 3,4,3,4\rangle$, $\langle 3,4,4,4\rangle$,
$\langle 3,3,3,3,4\rangle$, $\langle 4,6,8\rangle$,
$\langle 4,6,6\rangle$, $\langle 5,6,6\rangle$,
$\langle 3,10,10\rangle$,
$\langle 3,5,3,5\rangle$, $\langle 3,4,5,4\rangle$,
$\langle 3,3,3,3,5\rangle$, $\langle 4,6,10\rangle$.
In higher dimensions, the non-trivial semi-regular polytopes are
the Gosset family
$G_4,\ldots, G_8$ (subscripts indicate dimension)
and two 4-dimensional polytopes $P_{96}$ and $P_{720}$
(subscripts now indicate the number of vertices).
Coxeter ([C]) uses
$$\left\{{3 \atop 3,3}\right\}, s\{3,4,3\},
\left\{{3 \atop 3,5}\right\},
h\gamma_5, 2_{21}, 3_{21} \hbox{ and } 4_{21}$$
instead of our $G_4$, $P_{96}$, $P_{720}$,
$G_5$, $G_6$, $G_7$ and $G_8$, respectively.

The spectrum of a prism of $n$-gonal basis ([CDS]) is
$$\pm 1 + 2 \cos(2 \pi i/n),\quad i = 0, 1, \ldots, n-1$$
and the spectrum of an antiprism of $n$-gonal basis is
$$2\left(\cos(2\pi i/n) + \cos(\pi i/n)\right),
\quad i = 0, 1, \ldots, 2n-1,$$
as we will see in Section 2.  Below,
we list the characteristic polynomials of the non-trivial semi-regular
polytopes (i.e., the characteristic polynomials 
of the adjacency matrices of their graphs).
The polynomials are factored in $\ZZ[\tau]$, where
$\tau = (1+\sc)/2$ and $\taub = (1 -\sc)/2$;
the numerical values of non-integer eigenvalues are
listed in the same order as the factors, with semi-colons separating
roots of different factors.
By the Perron-Frobenius theorem ([Ga]),
the eigenvalue with largest module is simple and
equals the number of neighbours of a vertex,
since adjacency matrices have non-negative entries and are irreducible.

\bigskip

{
\parindent=0pt
\parskip=0pt
\raggedright
\tolerance=10000

\vbox{ $\langle 3,6,6\rangle$

$(X-3) (X-2)^3 X^2 (X+1)^3 (X+2)^3$ }

\bigskip

\vbox{ $\langle 3,4,3,4\rangle$

$(X-4) (X-2)^3 X^3 (X+2)^5$ }

\bigskip

\vbox{ $\langle 4,6,6\rangle$

$(X-3) (X-1)^3 (X+1)^3 (X+3)
(X^2-3)^2 (X^2 - 2X -1)^3 (X^2 + 2X -1)^3$

$(1.732051, -1.732051; 2.414214, -0.414214; -2.414214, 0.414214)$ }

\bigskip

\vbox{ $\langle 3,8,8\rangle$

$(X-3) (X-2)^3 (X-1) X^5 (X+1)^3 (X+2)^5 (X^2 - X - 4)^3$

$(2.561553, -1.561553)$ }

\bigskip

\vbox{ $\langle 3,4,4,4\rangle$

$(X-4) (X-3)^3 (X-1)^2 X^4 (X+1)^6 (X+3)^2 (X^2 + X - 4)^3$

$(1.561553, -2.561553)$ }

\bigskip

\vbox{ $\langle 3,3,3,3,4\rangle$

$(X- 5) (X + 1)^4 (X^2 + 2 X - 2)^2 (X^2 - 2 X - 6)^3
(X^3 + X^2 - 4 X - 2)^3$

$(0.7320508, -2.7320508; 3.645751, -1.645751;
1.81361, -0.470683, -2.34292)$ }

\bigskip

\vbox{ $\langle 4,6,8\rangle$

$(X-3) (X-2)^2 (X-1)^4 X^4 (X+1)^4 (X+2)^2 (X+3)
(X^2 - 2X - 2)^3 (X^2 + 2X -2)^3 \allowbreak
{(X^3 + X^2 -4 X -2)^3}, {(X^3 - X^2 -4 X +2)^3}$

$(2.732051, -0.732051; 0.732051, -2.732051; \allowbreak
1.81361, -0.470683, -2.34292; 2.34292, 0.470683, -1.81361)$}

\bigskip

\vbox{ $\langle 3,5,3,5\rangle$

$(X-4) (X-2)^5 (X-1)^4 (X+1)^4 (X+2)^{10}
(X - 2\tau)^3 (X - 2\taub)^3$

$(3.236068; -1.236068)$ }

\bigskip

\vbox{ $\langle 5,6,6\rangle$

$(X-3) (X-1)^9 (X+2)^4
(X + 2 - \tau)^3 (X + 2 - \taub)^3
(X + \tau)^5 (X + \taub)^5 \allowbreak
{(X^2 - (1 + \tau)X - 2 + \tau)^3}
(X^2 - (1 + \taub)X - 2 + \taub)^3
(X^2 + X - 4)^4 (X^2 - X - 3)^5$

$(-0.381966; -2.618034; -1.618034; 0.618034;
2.756598, -0.138564; 1.820249, -1.438283; \allowbreak
1.561553, -2.561553; 2.302776, -1.302776)$ }

\bigskip

\vbox{ $\langle 3,10,10\rangle$

$(X-3) X^{10} (X+2)^{11}
(X - \tau)^4 (X - \taub)^4
(X^2 - X -2 - 2\tau)^3 (X^2 - X - 2 - 2\taub)^3 \allowbreak
{(X^2 - X - 3)^4} {(X^2 - X - 4)^5}$

$(1.618034; -0.618034; 2.842236, -1.842236; 1.506942, -0.506942;
2.302776, -1.302776; \allowbreak
1.561553, -2.561553)$ }

\bigskip

\vbox{ $\langle 3,4,5,4\rangle$

$(X-4) (X-1)^4 X^6 (X+1)^4
(X - 2 - \tau)^3 (X - 2 - \taub)^3
(X + 2 - \tau)^8 (X + 2 - \taub)^8 \allowbreak
{(X + 1 - 2\tau)^4} {(X + 1 - 2\taub)^4}
(X^3 - X^2 - 7 X + 4)^5$

$(3.618034; 1.381966; -0.381966; -2.618034; 2.236068; -2.236068;
\allowbreak 2.92542, 0.551929, -2.47735)$ }

\bigskip

\vbox{ $\langle 3,3,3,3,5\rangle$

$(X-5) (X+1)^6
(X^2 - 2\tau X - 4 - \tau)^3 (X^2 -2\taub X - 4 - \taub)^3
(X^4 - 8 X^2 - 2 X + 10)^4 \allowbreak
{(X^5 + X^4 - 11 X^3 - 19 X^2 - X + 1)^5}$

$(4.48789, -1.25182; 1.32205, -2.55812;
2.71687, 1.07082, -1.50739, -2.2803; \allowbreak
3.5766, 0.195279, -0.285153, -2.1357, -2.35102)$ }

\bigskip

\vbox{ $\langle 4,6,10\rangle$

$(X-3) (X-1)^6 (X+1)^6 (X+3)
(X^2 + 2 X - 1 - \tau)^3 (X^2 - 2 X - 1 - \tau)^3 \allowbreak
{(X^2 + 2 X - 1 - \taub)^3 (X^2 - 2 X - 1 - \taub)^3} \allowbreak
{(X^4 - 6 X^2 - 2 X + 2)^4 (X^4 - 6 X^2 + 2 X + 2)^4} \allowbreak
{(X^5 - 3 X^4 - 3 X^3 + 11 X^2 - X - 3)^5
(X^5 + 3 X^4 - 3 X^3 - 11 X^2 - X + 3)^5}\allowbreak$

$(0.902113, -2.902113; 2.902113, -0.902113;
0.175571, -2.175571; 2.175571, -0.175571; \allowbreak
2.54501, 0.439406, -0.830209, -2.15421; \allowbreak
2.15421, 0.830209, -0.439406, -2.54501; \allowbreak
2.72142, 1.88838, 0.684645, -0.466437, -1.82801; \allowbreak
1.82801, 0.466437, -0.684645, 1.88838, 2.72142)$ }

\bigskip

\vbox{ $G_4$

$(X - 6) (X - 1)^4 (X + 2)^5$ }

\bigskip

\vbox{ $P_{96}$

$(X - 9) (X - 3)^8 (X - 1)^8 X^{14} (X + 2)^{24} (X + 3)^6
(X^2 - 4 X - 24)^4 (X^3 - X^2 - 16 X - 16)^9$

$(7.291503, -3.291503; 4.91638, -1.19656, -2.71982)$ }

\bigskip

\vbox{ $P_{720}$

${(X - 10) (X - 3)^{16} X^{16} (X + 1)^{60} (X + 2)^{28}} \allowbreak
{(X - \sc)^{24} (X + \sc)^{24}} \allowbreak
{(X - 5 - 2 \sc)^4 (X - 5 + 2 \sc)^4} \allowbreak
{(X + 1 - 3 \tau)^{24} (X + 1 - 3 \taub)^{24}} \allowbreak
{(X + 2 + \tau)^{16} (X + 2 + \taub)^{16}} \allowbreak
{(X + 1 + \tau)^{48} (X + 1 + \taub)^{48}} \allowbreak
{(X^2 + 3 X - 3)^{40}} \allowbreak
{(X^2 - 7 X - 4)^{16}} \allowbreak
{(X^2 + (-3 + 2\sc)X - 10)^9 (X^2 + (-3 - 2\sc)X - 10)^9} \allowbreak
{(X^2 - X -10 + 4\sc)^{36} (X^2 - X -10 - 4\sc)^{36}} \allowbreak
(X^3 - 4 X^2 - 15 X + 6)^{25}$

$(2.23607; -2.23607; 9.47214; 0.527864;
3.8541; -2.8541; -3.61803; -1.38197; \allowbreak
-2.61803; -0.38197; \allowbreak
0.79129, -3.79129; \allowbreak
7.53113, -0.53113; \allowbreak
2.51075, -3.98288; 8.63078, -1.15864; \allowbreak
1.642685, -0.642685; 4.88113, -3.88113; \allowbreak
6.2473, 0.36732, -2.61463)$ }

\bigskip

\vbox{ $G_5$

$(X - 10) (X - 2)^5 (X + 2)^{10}$ }

\bigskip

\vbox{ $G_6$

$(X - 16) (X - 4)^6 (X + 2)^{20}$ }

\bigskip

\vbox{ $G_7$

$(X - 27) (X - 9)^7 (X + 1)^{27} (X + 3)^{21}$ }

\bigskip

\vbox{ $G_8$

$(X - 56) (X - 28)^8 (X - 8)^{35} (X + 2)^{112} (X + 4)^{84}$ }

}

\bigbreak
{\nobigbf 2. Archimedean solids I}

\nobreak\bigskip\nobreak

In this section, we compute the spectra of seven
Archimedean solids:
$\langle 3,6,6\rangle$, $\langle 3,4,3,4\rangle$,
$\langle 4,6,6\rangle$, $\langle 3,8,8\rangle$,
$\langle 3,4,4,4\rangle$, $\langle 4,6,8\rangle$
and $\langle 3,5,3,5\rangle$ as well as the spectra of the antiprisms,
by using mirrors ([PS]) and relations between adjacency matrices.
We illustrate this method by computing
the spectrum of $\langle 3,4,4,4\rangle$.

\midinsert\smallskip

\centerline{\psfig{file=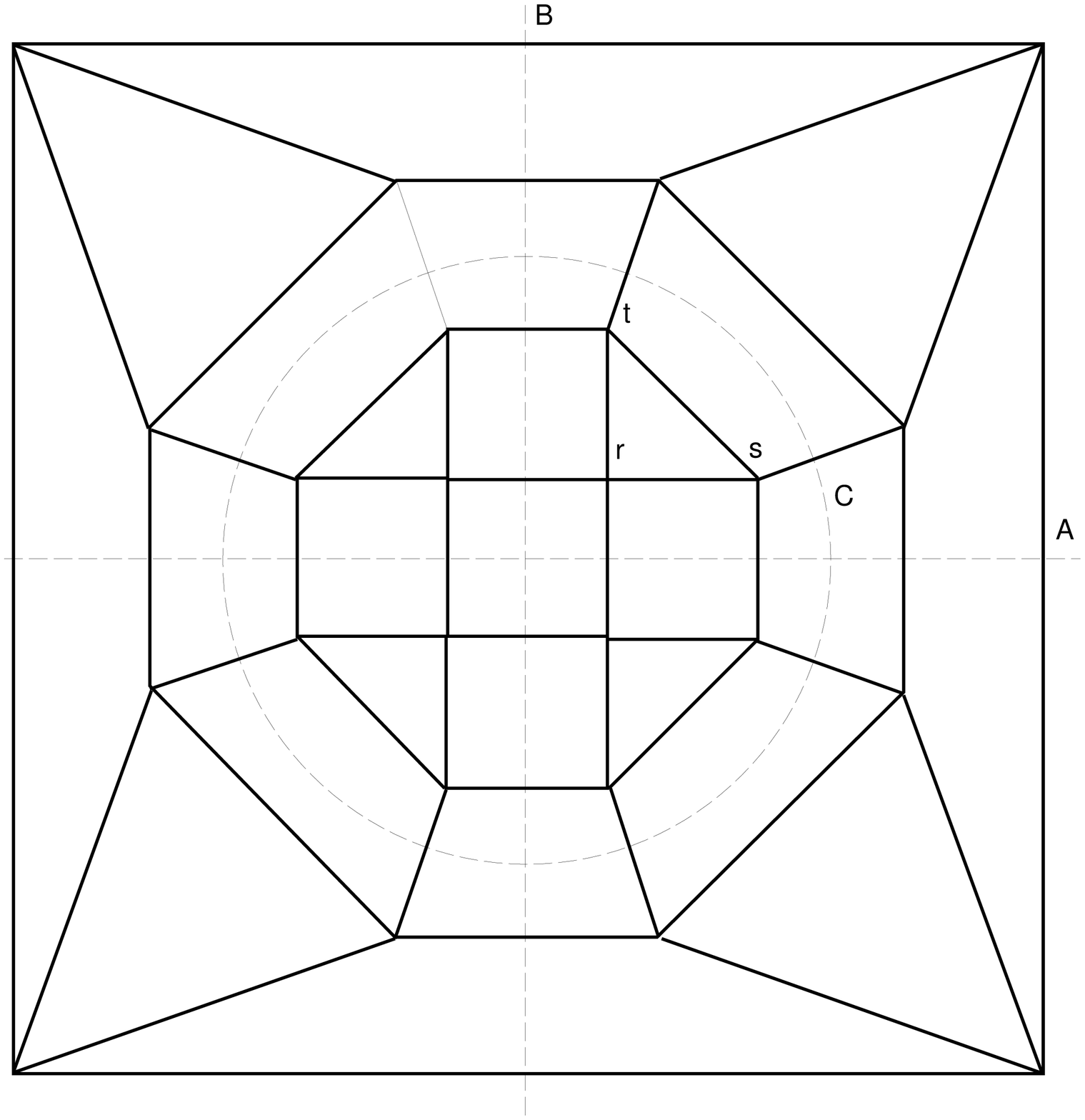,width=3.0in}}

\smallskip

\centerline{\eightbf Figure 2.1}

\smallskip\endinsert

In Figure 2.1,
full lines form the Schlegel diagram ([C])
of $\langle 3,4,4,4\rangle$, i.e.,
the stereographic projection of the polyhedron to the plane. Dotted
horizontal, vertical and round
lines indicate three planes of symmetry associated to three commuting
involutions $A, B$ and $C$.
Let $V$ be the vector space of complex valued functions
on the set of vertices of $\langle 3,4,4,4\rangle$.
Let $X$ be the $24 \times 24$ adjacency matrix of the polyhedron
(i.e., of its graph).
The involutions $A$, $B$ and $C$, as well as the matrix $X$,
can be interpreted as commuting linear transformations from $V$ to $V$.
Thus, the eigenspaces $V_+^A$ and $V_-^A$ associated to the eigenvalues
$1$ and $-1$ of $A$ are invariant under $B$, $C$ and $X$.
In a similar fashion, $V$ splits as a direct sum of eight subspaces
(which a priori could be trivial) of the form $V_{s_A,s_B,s_C} =
V_{s_A}^A \cap V_{s_B}^B \cap V_{s_C}^C$, where subscripts denote signs.
These subspaces are invariant under $X$ and the problem of computing the
spectrum of $X$ reduces to the same problem for the restriction
$X_{s_A,s_B,s_C}$ of $X$ to $V_{s_A,s_B,s_C}$.

Each invariant subspace can be coordinatized by the values $r$, $s$ and
$t$ of each vector on the three vertices indicated in Figure 2.1:
mirroring by $A$, $B$ and $C$ prescribes the values at the remaining
vertices. In particular, in this example, the eight invariant subspaces
are of dimension $3$. The adjacency matrices for
$V_{++-}, V_{+-+}$ and $V_{-++}$ are trivially conjugate (by renaming
the dotted lines), and therefore have the same spectrum; the same
happens to $V_{+--}, V_{-+-}$ and $V_{--+}$. We have
$$X_{+++} = \pmatrix{2 & 1 & 1 \cr 1 & 2 & 1 \cr 1 & 1 & 2 \cr}, \qquad
  X_{++-} = \pmatrix{2 & 1 & 1 \cr 1 & 0 & 1 \cr 1 & 1 & 0 \cr},$$
$$X_{+--} = \pmatrix{0 & 1 & 1 \cr 1 & 0 & 1 \cr 1 & 1 & -2 \cr}, \qquad
  X_{---} = \pmatrix{-2 & 1 & 1 \cr 1 & -2 & 1 \cr 1 & 1 & -2 \cr},$$
and the spectrum of $X$ is thus immediately computed.

For $\langle 3,6,6\rangle$, there are two commuting involutions
corresponding to reflections with respect to the two orthogonal planes
indicated by dotted lines in the Schlegel diagram in Figure 2.2.
The values $r$, $s$, $t$ and $u$ in the vertices marked in the figure
determine a unique vector in $V_{++}$.
Similarly, $r$, $s$, $t$ provide a basis for $V_{+-}$
and $s$, $t$ a basis for $V_{--}$
($V_{-+}$ is equivalent to $V_{+-}$ by renaming planes).

\midinsert\smallskip

\centerline{\psfig{file=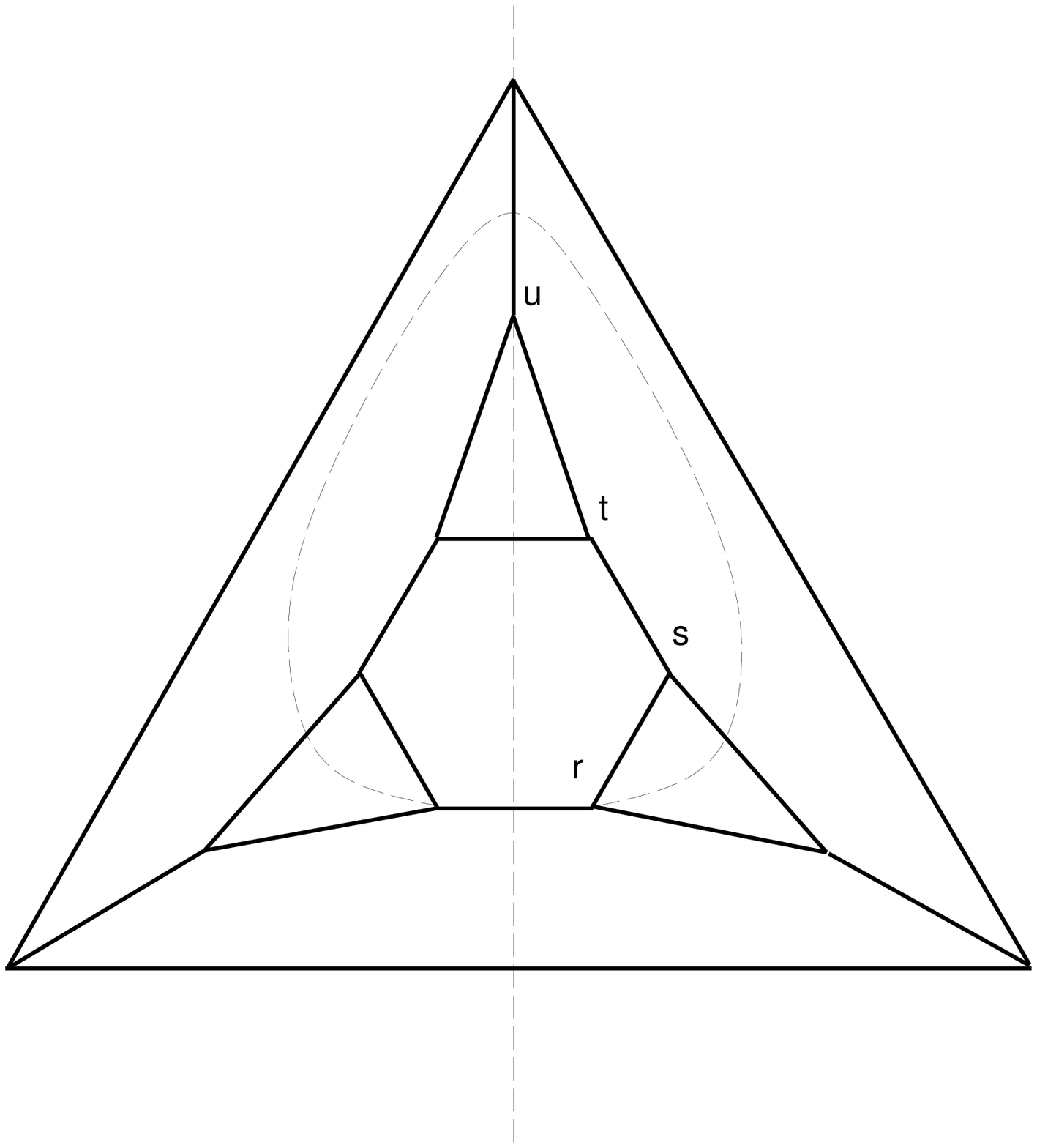,width=3.0in}}

\smallskip

\centerline{\eightbf Figure 2.2}

\smallskip\endinsert

The same three mirrors used for $\langle 3,4,4,4\rangle$
work for $\langle 3,4,3,4\rangle$,
reducing the problem to that of finding the spectra of a
$3 \times 3$ matrix in $V_{+++}$, a $2 \times 2$ matrix in $V_{++-}$
and a $1 \times 1$ matrix in $V_{+--}$
($V_{---}$ is 0-dimensional).
Again, the mirrors used for $\langle 3,4,4,4\rangle$
convert the computation
of the spectrum of $\langle 3,8,8\rangle$ into the study of
the restriction of the adjacency matrices to invariant subspaces
of dimension 3.
For $\langle 4,6,8\rangle$, however,
these mirrors produce $6 \times 6$ matrices.
A fourth mirror, indicated in Figure 2.3,
can be used to further split the invariant subspaces.
The reflection $D$ on the fourth mirror
does not commute with either $A$ or $B$.
Still, $V_{+++}, V_{++-}, V_{--+}$ and $V_{---}$
are invariant under $D$ and can thus be decomposed
into the two eigenspaces for the restrictions of $D$.
Notice that the four remaining subspaces are not invariant
under $D$. This is no serious problem, however:
for example, the restriction of the adjacency matrix to $V_{-++}$
has the same spectrum as the restriction to $V_{++-}$,
and this last space can be split by $D$.
The values $r, s, t$ in the figure again provide a basis
for each of the eight relevant invariant subspaces
and the spectrum of $\langle 4,6,8\rangle$ is now easily obtained.
The solid $\langle 4,6,6\rangle$ can be handled in the same fashion;
in section 5 we compute its spectrum using other methods.

\midinsert\smallskip

\centerline{\psfig{file=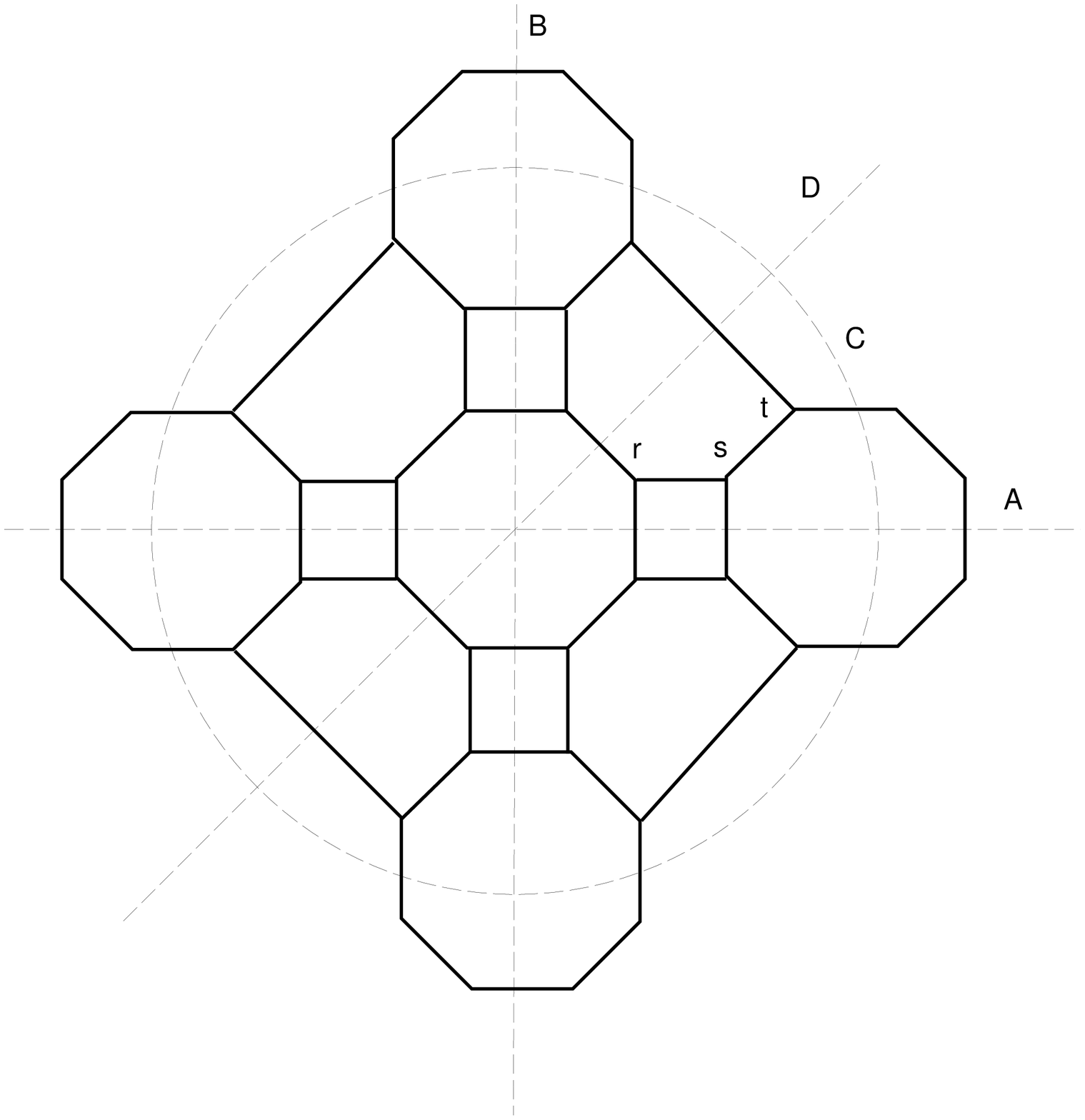,width=3.0in}}

\smallskip

\centerline{\eightbf Figure 2.3}

\smallskip\endinsert

We finish this section with three examples of computations of spectra
based on algebraic relations between adjacency matrices.
Let $Y$ be the adjacency matrix of a $2n$-gon:
then $X = Y^2 + Y - 2I$ is the adjacency matrix of the antiprism
with $n$-gonal basis 
and an eigenvalue $\lambda$ of $Y$ corresponds to an eigenvalue
$\lambda^2 + \lambda - 2$ of $X$.
As another example, the spectrum of $X_{3434}$,
the adjacency matrix of $\langle 3,4,3,4\rangle$,
can be computed in terms of the spectrum of $X_{444}$,
the adjacency matrix of the cube.
Indeed, let $Y_{v,e}$ be a $8 \times 12$ incidence matrix obtained
by numbering vertices and edges of the cube and setting
$y_{i,j} = 1$ if the $i$-th vertex belongs to the $j$-th edge.
Dually, let $Y_{e,v} = {(Y_{v,e})}^T$.
Notice that $X_{444} + 3 I = Y_{v,e} Y_{e,v}$
and that $X_{3434} + 2 I = Y_{e,v} Y_{v,e}$:
thus, the spectra of the left hand sides are equal
except for 4 extra zero eigenvalues in the second one.
The same method yields the spectrum of $\langle 3,5,3,5\rangle$
from the spectrum of the dodecahedron.
The spectrum of a regular polygon,
the cube and the dodecahedron are given in [CDS].

\bigbreak
{\nobigbf 3. Gosset polytopes}

\nobreak\bigskip\nobreak

In this section, we compute the spectra of 
the {\sl Gosset polytopes} $G_n$.
The solid $G_3$ is the triangular prism in dimension 3
and the vertex figure of $G_n$ is $G_{n-1}$.
Coordinates for the vertices of the Gosset polytopes 
are given in [G], [C] and [C3].

The polytope $G_n$ has for faces simplices $\alpha_{n-1}$
and cross polytopes $\beta_{n-1}$.
Recall that a cross polytope is a generalized octahedron,
or $\{3,\ldots,3,4\}$ in Schl\"afli's notation;
the vertices of $\beta_n$ consist of a north and south pole
together with the vertices of an equatorial $\beta_{n-1}$.
We make extensive use of Table 3.1 ([C], section 11.8).

\midinsert\smallskip

\centerline{\vbox{\offinterlineskip\halign{%
\strut\quad\hfil#\hfil\quad&\vrule#&&\strut\quad\hfil#\hfil\quad\cr
&depth 8pt&$G_4$&$G_5$&$G_6$&$G_7$&$G_8$\cr
\noalign{\hrule}
Number of vertices&height 12pt&10&16&27&56&240\cr
Number of $\beta_{n-1}$ faces&&5&10&27&126&2160\cr}}}

\smallskip

\centerline{\bf Table 3.1}

\smallskip\endinsert

A more combinatorial description of the graphs
$G_i$, $i = 4,\ldots,8$, is given in [BCN]
where the notation $E_i(1)$
(reminiscent of the Lie algebras $E_n$) is employed.

Let $\Gamma(G_n)$ be the group of isometries of $G_n$
and, for an arbitrary choice of a vertex $p_1$ of $G_n$,
call $\Gamma_1(G_n)$ the subgroup of $\Gamma(G_n)$ fixing $p_1$.
As explained in [C3], $\Gamma_1(G_n) = \Gamma(G_{n-1})$.
Also, $\Gamma(G_n)$ acts transitively both on
the set of simplicial faces and
on the set of cross polytope faces.
The orbits of the action of $\Gamma_1(G_n)$
on the set of vertices of $G_n$ are called
{\sl suborbits} in the terminology of association schemes
and {\sl levels} in [ST].
The vertex $p_1$ forms a suborbit by itself;
the next suborbit is a $G_{n-1}$.
As we shall see, the number of suborbits of $G_n$
is 3, 3, 3, 4 and 5 for $n = 4$, 5, 6, 7 and 8.
Thus, the method employed in [ST] for regular polytopes
is well suited for the Gosset family.
Actually, this method is a special case of the technique
of {\sl regular partitions} ([BCN]) for the computation
of the spectrum of an association scheme.

We briefly describe this special method.
Let $v \ne 0$ satisfy $X v = \lambda v$
for the adjacency matrix $X$ of a semi-regular polytope.
Without loss, the value of $v$ at $p_1$ is nonzero.
Let $S$ be the space of vectors which are constant on suborbits:
$S$ is invariant under $X$.
The (nonzero) average $s \in S$ of $v$ under the action of $\Gamma_1$
is another eigenvector of $X$ associated
to the same eigenvalue $\lambda$.
Thus, up to multiplicities, the spectrum of $X$ equals
the spectrum of its restriction $B$ to $S$.
For the basis of $S$ consisting 
of vectors $s_k$ taking the value 1 at the $k$-th suborbit
and 0 elsewhere, the entry $b_{ij}$ of $B$
is the number of neighbours in suborbit $j$ of a vertex in suborbit $i$.
We now obtain the matrices $B_n$ for $G_n$, $n \ge 4$.

The 16 vertices of $G_5$ are $(\pm 1,\ldots,\pm 1)$
with an even number of minus signs ([C]).
Taking $p_1 = (1,\ldots,1)$, different sums of coordinates
clearly imply different suborbits.
The 3 possible values of this sum, 5, 1 and $-3$,
correspond to sets of 1, 10 and 5 vertices
which might, in principle, be the union of more than one suborbit.
The intermediate set is a $G_4$,
on which $\Gamma_1(G_5) = \Gamma(G_4)$ is transitive:
this set is therefore a suborbit.
Each of the five $\beta_3$ faces of the second suborbit
is an equator of a $\beta_4$ face of $G_5$
with north pole $p_1$ and south pole in the third set.
Since $\Gamma(G_4)$ acts transitively
on the cross polytope faces of $G_4$, this bottom set is a suborbit.
The $3 \times 3$ matrix $B_5$ is
$$B_5 = \pmatrix{0 & 10 & 0 \cr 1 & 6 & 3 \cr 0 & 6 & 4 \cr}$$
with eigenvalues 10, 2 and $-2$.
Let their multiplicities in $X$ be 1, $a$ and $b$
(recall that the largest eigenvalue is always simple).
We must then have $1 + a + b = 16$ and
$1 \cdot 10 + a \cdot 2 + b \cdot (-2) = \tr(X) = 0$;
$a$ and $b$ are then 5 and 10.

The central suborbit of $G_5$ gives coordinates (in $\RR^5$)
for $G_4$ and, as above, it is easy to obtain
the decomposition of $G_4$ into suborbits of 1, 6 and 3 elements,
yielding
$$B_4 = \pmatrix{0 & 6 & 0 \cr 1 & 3 & 2 \cr 0 & 4 & 2 \cr},$$
from which the spectrum of $G_4$ follows.

From coordinates for $G_5$,
one obtains coordinates for a $G_6$ inscribed in $S^6$
for which suborbits are indicated by the first entry.
Start with $p_1 = (1,0,\ldots,0)$ and add the next suborbit of vertices,
of the form $(a, b\,G_5)$,
where $a = 1/4$ and $b = (\sqrt{3}/4)$;
the coefficients are taken so that the second suborbit
lies in the unit sphere and all edges have equal length.
The vertices of a $\beta_5$ face of $G_6$ containing $p_1$
are, besides $p_1$, the vertices of a $\beta_4$ face
of the second suborbit (a $G_5$) and the reflection of $p_1$ on
the hyperplane containing this $\beta_4$.
We thus obtain 10 new vertices of $G_6$, with first coordinate $-1/2$,
associated to the 10 $\beta_4$ faces of the second suborbit;
$\Gamma_1(G_6)$ is transitive on the set of $\beta_4$ faces
of the second suborbit and these 10 points form a third suborbit of $G_6$.
Since $G_6$ has 27 vertices, these are the only suborbits;
adjacencies among suborbits are given by
$$B_6 = \pmatrix{0 & 16 & 0 \cr 1 & 10 & 5 \cr 0 & 8 & 8 \cr}$$
and we are done with $G_6$.

In a similar fashion, coordinates for $G_7$
in which the first entry indicates the suborbit are
$$(1,0,\ldots,0),\;(1/3,(2\sqrt{2}/3)\,G_6),\;%
(-1/3,-(2\sqrt{2}/3)\,G_6),\;(-1,0,\ldots,0).$$
The polytope $G_7$ then admits an involution taking $v$ to $-v$
splitting $V_{G_7}$ into two 28-dimensional subspaces
$V_+$ and $V_-$, invariant under $X$, and
$S$ into 2-dimensional subspaces $S_+$ and $S_-$, invariant under $B_7$.
The two restrictions of $B_7$ are
$$B_+ = \pmatrix{0 & 27 \cr 1 & 26 \cr},\quad
B_- = \pmatrix{0 & 27 \cr 1 & 6 \cr},$$
with spectra 27, $-1$ and 9, $-3$.
Considered as eigenvalues of $X$,
27 and $-1$ have multiplicities 1 and 27
since their eigenspaces are contained in $V_+$.
The two remaining multiplicities are now computed
taking into account that $\tr(X) = 0$.

Finally, analogous coordinates for $G_8$ are
$$(1,0,\ldots,0),\;(1/2,(\sqrt{3}/2)\,G_7),\;%
(0,\ast),\;(-1/2,(\sqrt{3}/2)\,G_7),\;(-1,0,\ldots,0),$$
where the vertices of the central suborbit
correspond to the 126 $\beta_6$ faces of the second suborbit (a $G_7$).
Again, $G_8$ admits an antipodal involution and $B_8$ splits as
$$B_+ = \pmatrix{0 & 56 & 0 \cr 1 & 28 & 27 \cr 0 & 24 & 32 \cr},\quad
B_- = \pmatrix{0 & 56 \cr 1 & 26 \cr}.$$
To compute the multiplicities, we only need the additional remark
that, since antipodal points are never neighbours,
the restrictions of $X$ to the 120-dimensional spaces $V_+$ and $V_-$
have trace zero.

\bigbreak
{\nobigbf 4. Groups, Cayley graphs and representations}

\nobreak\bigskip\nobreak

By using representations of groups,
Lov\'asz ([L]) described a method to compute
the spectrum of a Cayley graph, or, more generally,
of a graph with a transitive group of isomorphisms.
Babai ([B]) modified Lov\'asz's method
(and corrected a minor mistake in [L]) and
obtained families of Cayley graphs with the same spectrum.
We present yet another technique based on representations of groups
which is more convenient for our purposes.

For an undirected graph $\sG$, the {\sl doubling} $\tilde\sG$ of $\sG$
is the directed graph obtained by substituting
each edge of $\sG$ by two edges with opposite orientations.
A {\sl Cayley structure} for an undirected graph $\sG$
consists of a choice of a vertex (to be the identity)
and a colouring of the edges of the doubling $\tilde\sG$ of $\sG$
(to represent the set $H$ of generators)
such that $\tilde\sG$ becomes a Cayley graph, i.e.,
the group of isomorphisms of the coloured graph
is simply transitive on vertices.
Notice that $h \in H$ implies $h^{-1} \in H$.
Similarly, a {\sl Cayley structure} for a polytope
is a Cayley structure for its 1-skeleton.
Most semi-regular polytopes admit Cayley structures: in dimension 3,
only the dodecahedron and $\langle 3,5,3,5\rangle$ do not ([C1]).

Let $V_{\sG}$ be the complex vector space of functions
from the vertices of the Cayley graph $\sG$ to $\CC$.
Define $e_k$ by $e_k(g) = \delta_{kg}$.
The {\sl canonical left representation} of $\Gamma$ on $V_{\sG}$
is given by $L_g e_k = e_{gk}$ or,
equivalently, $(L_h v)(g) = v(h^{-1}g)$.
For a polytope with a Cayley structure,
the canonical left action corresponds
to an action by isometries.
The {\sl canonical right representation} is
given by $R_g e_k = e_{kg^{-1}}$.
The two canonical representations are equivalent, being
intertwined by the linear involution $e_k \mapsto e_{k^{-1}}$,
and commute: $L_{g_1} R_{g_2} = R_{g_2} L_{g_1}$.
Geometrically, this action can be interpreted in terms of
the Cayley graph $\sG$ as follows: 
for a generator $h$, the value of $R_h v$ at a vertex $g$
is the value of $v$ at the target vertex of the $h$-edge
starting from $g$ (which is, of course, $gh$),
i.e., $(R_h v)(g) = v(gh)$.
Thus, unless $\Gamma$ is abelian,
$R_h$ does not preserve adjacencies among vertices.

From representation theory for finite groups ([Se]),
the space $V_{\sG}$ splits as a direct sum,
$$V_{\sG} =
\bigoplus_{1 \le i \le t, 1 \le j \le d_i} W_{i,j},$$
where each $W_{i,j}$ is invariant under the canonical right action
and the restriction of the action to $W_{i,j}$
is isomorphic to the $i$-th irreducible representation.

As in Babai ([B]),
the adjacency matrix $X$ of the oriented graph $\sG$,
defined by $x_{ij} = 1$ if there is an (oriented) edge
from the $i$-th to the $j$-th vertex, and 0 otherwise, satisfies
$$X = \sum_{h \in H}R_h.$$
In particular, $X$ commutes with the left canonical action.
Invariant subspaces for the canonical right representation
are therefore invariant under $X$ and the restriction $X_{i,j}$ of
$X$ to each $W_{i,j}$ is the sum
of the restrictions of the linear transformations
corresponding to the elements of $H$.
Thus, if we have matrices $M_{h,i}$ for the generators in $H$
in any representation isomorphic to
the $i$-th irreducible representation of $\Gamma$,
$X_{i,j}$ is conjugate to
$$X_i = \sum_{h \in H}M_{h,i}.$$
Clearly, the spectrum $\sigma(X)$ of $X$ is the union (as a multiset)
of $d_i$ copies of $\sigma(X_i)$, $i = 1, \ldots, t$.

Although we use the right canonical action
to decompose $V_{\sG}$, the left canonical action
is geometrically more helpful in the identification
of a Cayley structure for a polytope.
Like Babai's method,
this technique applies to polytopes with different weights
for edges of different colours.

Instead of computing $X_i$,
Lov\'asz and Babai have a formula for traces of powers of $X_i$ ([B]):
$$\tr(X_i^\ell) =
\lambda_{i,1}^\ell + \cdots + \lambda_{i,{d_i}}^\ell =
\sum_{h_1,\ldots,h_\ell \in H}\chi_i(h_1\cdot\;\cdots\;\cdot h_\ell),$$
where $\lambda_{i,1},\ldots,\lambda_{i,{d_i}}$
are the eigenvalues of $X_i$, with dimensions $d_i$,
and $\chi_i$, $i = 1, \ldots, t$,
are the characters of the irreducible representations of $\Gamma$.
The traces of $X_i^\ell$ for
$\ell = 1, \ldots, d_i$, determine $\sigma(X_i)$.

In our examples, we find the explicit computation
of the right hand side of the above formula cumbersome
since it involves counting special paths
and we prefer to specify matrices for each irreducible representation.

\bigbreak
{\nobigbf 5. Archimedean solids II}

\nobreak\bigskip\nobreak

In this section we compute the spectra of the six remaining
Archimedean solids using the methods of the previous section.

As a first example, let
$P = \langle 4,6,6\rangle$.
The cube and $P$ have the same group of symmetries (of order 48);
the subgroup $\Gamma$ of orientation preserving isometries
acts simply transitively on the vertices of $P$,
thus giving $P$ a Cayley structure.
As is well known, $\Gamma$ is isomorphic to $S_4$,
the permutations of $\{1,2,3,4\}$
(number pairs of antipodal hexagons $1, 2, 3, 4$ as in Figure 5.1).
With this notation for $\Gamma$, $H = \{ (12), (1234), (1432) \}$:
the first generator corresponds in Figure 5.1
to the edge joining $e$ to $a$,
the second to the edge going southwest from $e$ to $d$
and the third to the remaining edge starting from $e$.
Generators for the five irreducible representations
of $S_4$ are given in Table 5.2.
The corresponding $X_i = M_{(12),i} + M_{(1234),i} + M_{(1432),i}$,
with their spectra, are listed in Table 5.3.

\midinsert\smallskip

\centerline{\psfig{file=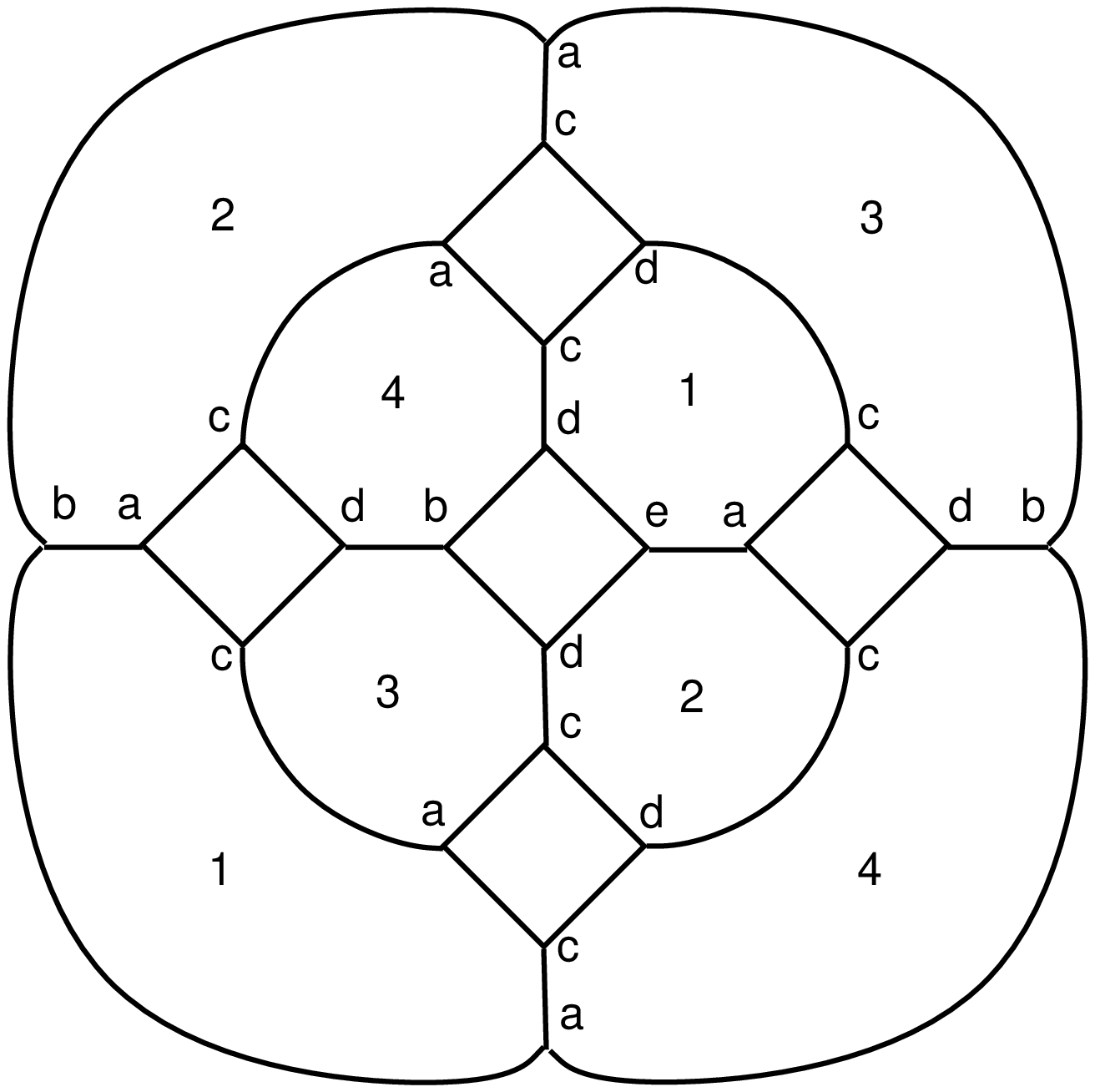,width=3.0in}}

\smallskip

\centerline{\eightbf Figure 5.1}

\smallskip\endinsert

\midinsert\smallskip

$$\matrix{
M_{(12),1} = 1& M_{(1234),1} = 1\cr
&\cr
M_{(12),2} = -1& M_{(1234),2} = -1\cr
&\cr
M_{(12),3} = \pmatrix{0 & 1 \cr 1 & 0 \cr} &
M_{(1234),3} = \pmatrix{-1 & 0 \cr -1 & 1 \cr} \cr
&\cr
M_{(12),4} =
\pmatrix{0 & 1 & 0 \cr 1 & 0 & 0 \cr 0 & 0 & 1 \cr}&
M_{(1234),4} =
\pmatrix{0 & 0 & -1 \cr 1 & 0 & -1 \cr 0 & 1 & -1 \cr}\cr
&\cr
M_{(12),5} = - M_{(12),4}&
M_{(1234),5} = - M_{(1234),4}\cr
}$$

\smallskip

\centerline{\eightbf Table 5.2}

\smallskip\endinsert

\midinsert\smallskip

$$\matrix{
X_1 = 3& \sigma(X_1) = \{3\}\cr
&\cr
X_2 = -3& \sigma(X_2) = \{-3\}\cr
&\cr
X_3 = \pmatrix{-2 & 1 \cr -1 & 2 \cr}&
\sigma(X_3) = \{\pm \sqrt{3}\}\cr
&\cr
X_4 = \pmatrix{-1 & 2 & -1 \cr 1 & 0 & 0 \cr -1 & 1 & 0 \cr}&
\sigma(X_4) = \{1, -1 \pm \sqrt{2}\}\cr
&\cr
X_5 = - X_4&
\sigma(X_5) = \{-1, 1 \pm \sqrt{2}\}\cr
}$$

\smallskip

\centerline{\eightbf Table 5.3}

\smallskip\endinsert

Consider now $\langle 3,3,3,3,4\rangle$,
which admits a Cayley structure
for its full symmetry group $\Gamma = S_4$.
Comparing Figures 5.1 and 5.4, we see two additional generators
$(234)$ and $(243)$, so
$H = \{ (12), (1234), (1432), (234), (243) \}$.
From Table 5.2 we have
$M_{(12),i}$ and $M_{(1234),i}$ and thus we have also
$M_{(234),i} = M_{(12),i}M_{(1234),i}$ and
$M_{(243),i} = M_{(234),i}^{-1}$.
The restrictions $X_i$ and their spectra follow immediately.

\midinsert\smallskip

\centerline{\psfig{file=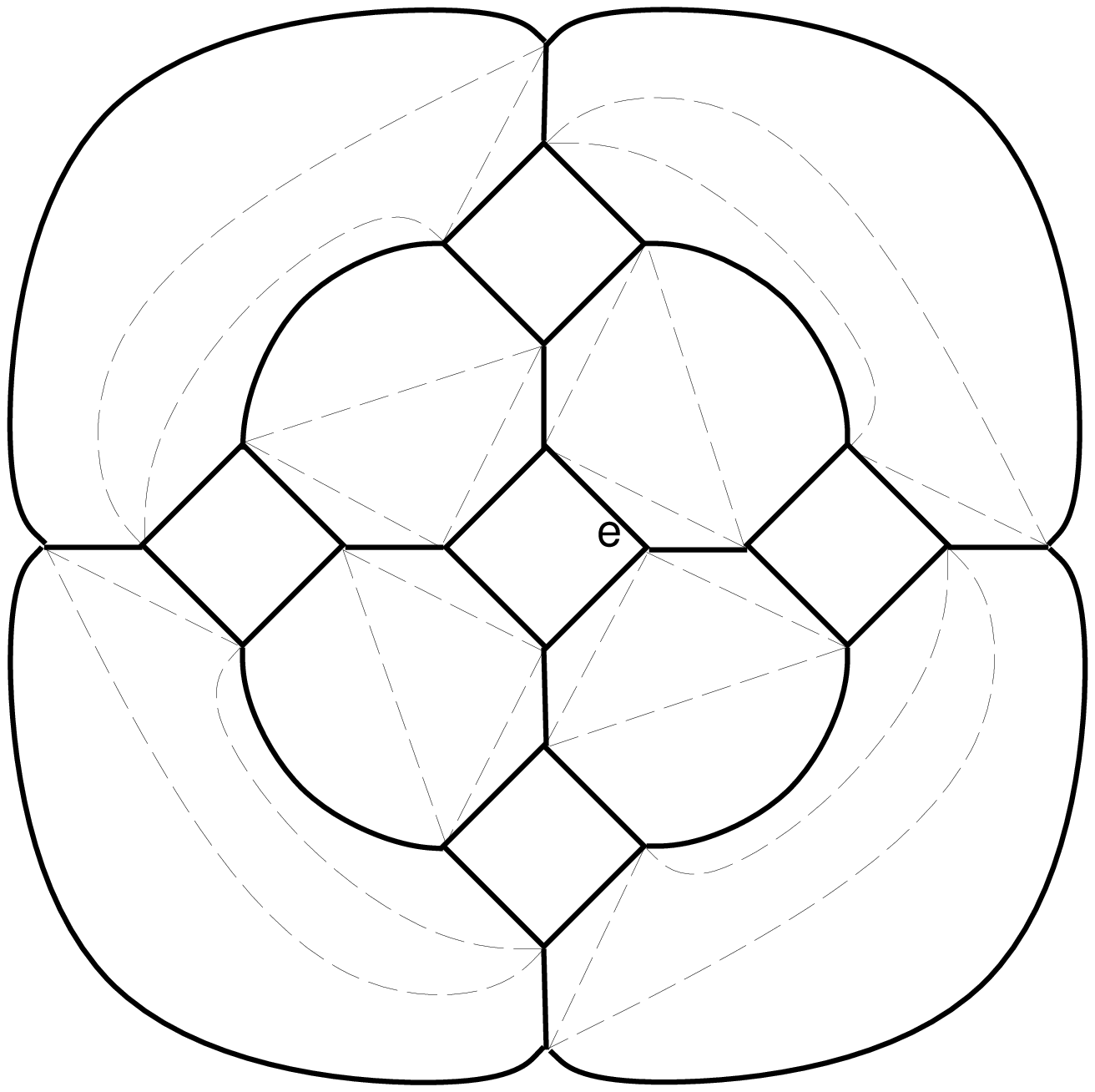,width=3.0in}}

\smallskip

\centerline{\eightbf Figure 5.4}

\smallskip\endinsert

To handle $\langle 3,4,5,4\rangle$,
consider the planes containing the triangular faces:
they form the faces of a circumscribed icosahedron.
The group $\Gamma$ of orientation preserving isometries
of these two polyhedra is therefore the same.
It is well known that this group can be identified
with $A_5$, the even permutations of $\{1,2,3,4,5\}$.
The polytope $\langle 3,4,5,4\rangle$ admits a Cayley structure
since $A_5$ acts simply transitively on its 60 vertices.
With an appropriate identification of $A_5$ with
the group of isometries of the icosahedron,
$H = \{a, a^{-1}, b, b^{-1}\}$, for $a = (12345)$ and $b = (253)$.
Generators for the five irreducible representations of $A_5$
are given in Table 5.5 where,
again, $\tau = (1 + \sqrt{5})/2$ and $\taub = (1 - \sqrt{5})/2$.
The matrices $X_i$ and their spectra are easily obtained.

\midinsert\smallskip

$$\matrix{
M_{a,1} = 1& M_{b,1} = 1\cr
&\cr
M_{a,2} = {1 \over 2}
\pmatrix{1 & \tau - 1 & -\tau \cr
\tau - 1 & \tau & 1 \cr \tau & -1 & \tau - 1 \cr}&
M_{b,2} = {1 \over 2}
\pmatrix{-1 & \tau - 1 & -\tau \cr
- \tau + 1 & \tau & 1 \cr \tau & 1 & - \tau + 1 \cr}\cr
&\cr
M_{a,3} = {1 \over 2}
\pmatrix{1 & \taub - 1 & -\taub \cr
\taub - 1 & \taub & 1 \cr \taub & -1 & \taub - 1 \cr}&
M_{b,3} = {1 \over 2}
\pmatrix{-1 & \taub - 1 & -\taub \cr
- \taub + 1 & \taub & 1 \cr \taub & 1 & - \taub + 1 \cr}\cr
&\cr
M_{a,4} =
\pmatrix{0 & 0 & 0 & -1 \cr 1 & 0 & 0 & -1 \cr
0 & 1 & 0 & -1 \cr 0 & 0 & 1 & -1 \cr}&
M_{b,4} =
\pmatrix{1 & -1 & 0 & 0 \cr 0 & -1 & 1 & 0 \cr
0 & -1 & 0 & 0 \cr 0 & -1 & 0 & 1 \cr}\cr
&\cr
M_{a,5} =
\pmatrix{0 & 0 & 0 & 0 & 1 \cr 1 & 0 & 0 & 0 & 0 \cr
0 & 1 & 0 & 0 & 0 \cr 0 & 0 & 1 & 0 & 0 \cr 0 & 0 & 0 & 1 & 0}&
M_{b,5} =
\pmatrix{-1 & 0 & 0 & 0 & 1 \cr -1 & 0 & 1 & 0 & 0 \cr
-1 & 0 & 0 & 1 & 0 \cr -1 & 1 & 0 & 0 & 0 \cr -1 & 0 & 0 & 0 & 0 \cr}\cr
}$$

\smallskip

\centerline{\eightbf Table 5.5}

\smallskip\endinsert

The three solids $\langle 5,6,6\rangle$,
$\langle 3,10,10\rangle$ and $\langle 3,3,3,3,5\rangle$
admit Cayley structures for the same $\Gamma = A_5$
with $H$ in each case being
$\{a, a^{-1}, ab\}$, $\{b, b^{-1}, ab\}$ and
$\{a, a^{-1}, b, b^{-1}, ab\}$.
The rest of the procedure is analogous.

The polyhedron $\langle 4,6,10\rangle$,
unlike the previous four examples,
has 120 vertices and admits a Cayley structure for its {\sl full}
isometry group, which is isomorphic to $A_5 \oplus \ZZ/(2)$,
where the generator of $\ZZ/(2)$ is the (orientation reversing)
antipodal map. The characters and representations for this
larger group are trivially obtained from those of $A_5$ ([Se])
and the ten $X_i$ matrices are thus easily obtained.
We prefer, however, to consider $X^2$,
the square of the adjacency matrix.
Since $\langle 4,6,10\rangle$ is bipartite,
$X^2$ has two obvious invariant
subspaces corresponding to the white and black vertices.
The white vertices are naturally identified
with those of $\langle 3,3,3,3,5\rangle$
and admit a simply transitive action of $A_5$.
Let $Y$ be the restriction of $X^2$ to the white vertices:
in the notation of section 4,
$Y = 3 I + R_a + R_{a^{-1}} + R_b + R_{b^{-1}} + 2 R_{ab}$.
Thus, the computation of $\sigma(Y)$ reduces as usual
to the computations of spectra of matrices in each irreducible
representation of $A_5$.
Finally, $\sigma(X)$ contains the two square roots
of each element of $\sigma(Y)$.

\bigbreak
{\nobigbf 6. Three discrete subgroups of $S^3$}

\nobreak\bigskip\nobreak

Among the six regular polytopes in four dimensions,
three of them, $\{3,3,4\}$, $\{3,4,3\}$ and $\{3,3,5\}$,
have a most remarkable property ([C], [C2]):
when suitably inscribed in the unit quaternionic sphere $S^3$,
their vertices form finite groups.
In this section, we briefly describe these three groups,
to be called $Q_8$, $Q_{24}$ and $Q_{120}$.

Following [C], we display appropriate choices
of unit quaternions for the vertices of these polytopes.
We identify $(a,b,c,d)$ to $a + bi + cj + dk$.
The elements of $Q_8$ are the vertices of $\{3,3,4\}$
with coordinates $\pm 1$, $\pm i$, $\pm j$ and $\pm k$.
In order to obtain the vertices of $\{3,4,3\}$,
add to the above list all 16 points of the form
$(\pm 1 \pm i \pm j \pm k)/2$, producing $Q_{24}$.
Finally, consider the 96 points obtained by even permutations of
the coordinates of $(\pm\tau,\pm 1,\pm\taub,0)/2$:
adding these to the 24 points already defined,
we obtain the vertices of $\{3,3,5\}$, i.e.,
the elements of $Q_{120}$.

The four groups $S^3$, $Q_8$, $Q_{24}$ and $Q_{120}$
have the same centre $\{\pm 1\}$.
As is well known ([MT]),
the quotient $S^3/\{\pm 1\}$ is isomorphic to $SO(3)$.
The quotient $Q_8/\{\pm 1\} = \ZZ/(2) \times \ZZ/(2)$,
as a group of isometries in $\RR^3$,
is generated by $180^\circ$ rotations around the three axis.
Similarly, $Q_{24}/\{\pm 1\} = A_4$ and
$Q_{120}/\{\pm 1\} = A_5$ are the groups of orientation preserving
isometries of a tetrahedron and an icosahedron, respectively.

Conjugacy classes in $S^3$ are determined by the first coordinate,
i.e., the real part;
in the groups $Q_n$, therefore,
points with different real parts are never conjugate.
The conjugacy classes of $Q_8$ are
$\{1\}$, $\{\pm i\}$, $\{\pm j\}$, $\{\pm k\}$, $\{-1\}$.
Representatives for the conjugacy classes of $Q_{24}$ are
1, $(1+i+j+k)/2$, $(1+i+j-k)/2$, $i$,
$(-1+i+j+k)/2$, $(-1+i+j-k)/2$ and $-1$;
these classes have 1, 8, 8, 6, 8, 8 and 1 elements, respectively.
For $Q_{120}$, the real part determines the conjugacy class:
there are thus 9 conjugacy classes denoted by
$1$, $b$, $c$, $d$, $e$, $-d$, $-c$, $-b$, $-1$ (decreasing real parts)
with 1, 12, 20, 12, 30, 12, 20, 12, 1 elements, respectively;
the order of an element in each conjugacy class is
1, 10, 6, 5, 4, 10, 3, 5, 2.
The conjugacy class $b$ consists of the 12 neighbours
of the vertex 1 and these 12 points form the vertices
of an icosahedron.

The irreducible representations of $S^3$ are well known ([Su]):
we call them $\sR_n$, where $\sR_n$ has dimension $n+1$.
The representation $\sR_0$ is trivial and
$\sR_1$ corresponds to the identification $S^3 = SU(2)$.
A basis for the space where $\sR_n$ acts is given by monomials
of degree $n$ in two variables $e_1$ and $e_2$,
identified with the basis of $\CC^2$;
for a homogeneous polynomial $\phi$ of degree $n$, define
$\sR_n(g)(\phi(e1,e2)) = \phi(\sR_1(g)(e_1),\sR_1(g)(e_2))$.
Restrictions of $\sR_n$ to $Q_n$ are still representations,
but are usually not irreducible.
Call an irreducible representation of $S^3$ or $Q_n$
{\sl even} (resp., {\sl odd}) if $-1$ is taken to $I$ (resp., $-I$);
$\sR_n$ is even if and only if $n$ is even.
The number of even (resp. odd) irreducible representations of
$Q_8$, $Q_{24}$ and $Q_{120}$ are 4, 4, 5 (resp. 1, 3, 4);
character tables for $Q_8$ and $Q_{24}$ are easy to obtain
and that of $Q_{120}$ is in [CCNPW].

We show how to obtain a matrix form for $\sM_i(G)$,
the $i$-th irreducible representation of $G$
(even representations come first;
within each parity, representations are ordered by dimension).
By restrictions, $\sR_1$ yield
$\sM_5(Q_8)$, $\sM_5(Q_{24})$ and $\sM_6(Q_{120})$.
Also, $\sR_2$ yields $\sM_4(Q_{24})$ and $\sM_2(Q_{120})$;
$\sR_3$, $\sR_4$ and $\sR_5$ yield
$\sM_8(Q_{120})$, $\sM_5(Q_{120})$ and $\sM_9(Q_{120})$.
The other representations can be obtained by algebraic conjugation
and tensor products:
$\sM_3(Q_{120})$ and $\sM_7(Q_{120})$ are the conjugates of
$\sM_2(Q_{120})$ and $\sM_6(Q_{120})$, respectively.
Also, $\sM_6(Q_{24}) = \sM_2(Q_{24}) \otimes \sM_5(Q_{24})$,
$\sM_7(Q_{24}) = \sM_3(Q_{24}) \otimes \sM_5(Q_{24})$ and
$\sM_4(Q_{120}) = \sM_6(Q_{120}) \otimes \sM_7(Q_{120})$.

In our examples, we frequently consider isometry groups in $\RR^4$.
As is well known ([MT]), $SO(4) = (S^3 \times S^3)/(-1,-1)$
by unit quaternion bilateral multiplication:
$(q,r) \cdot v = qvr^{-1}$.
Thus, the finite groups $(Q_n \times Q_m)/(-1,-1)$
act on $\RR^4$ by isometries.
The irreducible representations of such groups
are obtained by tensoring representations of each factor
having the same parity.

\bigbreak
{\nobigbf 7. $P_{96}$}

\nobreak\bigskip\nobreak

In this section, we compute the spectrum of $P_{96}$
(Coxeter's $s\{3,4,3\}$),
one of the three semi-regular polytopes of dimension 4.
Each of its 96 vertices is surrounded by
three icosahedra and five tetrahedra,
arranged according to the vertex figure shown in Figure 8.9A of [C];
each vertex has nine neighbours.
Our main task is to obtain a group $\Gamma$ of isometries of $P_{96}$
acting simply transitively on vertices.

We remind the reader of a construction for $P_{96}$,
detailed in [C], section 8.4.
Edges of $\{3,4,3\}$ can be oriented so that,
given any vertex $p$, there are four edges pointing outwards from $p$,
no two of which belong to the same (triangular) 2-cell.
Now, divide each oriented edge in two segments $a$ and $b$
(in this order) satisfying $b/a = \tau$;
the points thus obtained are the vertices of a $P_{96}$.
Alternatively, the 120 vertices of $\{3,3,5\}$
are the disjoint union of the 24 vertices of a $\{3,4,3\}$
and the 96 vertices of a $P_{96}$;
from the coordinates for $\{3,3,5\}$ in the previous section,
we thus obtain, as in [C] (section 8.7) coordinates for $P_{96}$.

The finite group $(Q_{24} \times Q_{24})/(-1,-1)$
acts on $\{3,4,3\}$ by isometries preserving edge orientation,
therefore acting also on $P_{96}$.
This group is too large to act simply transitively on the vertices
of $P_{96}$:
the subgroup $\Gamma = (Q_{24} \times Q_8)/(-1,-1)$ has the right order.
Clearly, $\Gamma$ acts transitively on the vertices of $\{3,4,3\}$
by elements of the form $(q,1)$.
Also, $\Gamma$ acts transitively on oriented edges:
by the transitivity on vertices, it is enough to show that it acts
transitively on the four edges starting from $1$,
which is done by elements of the form $(r,r)$.
Adding up, $\Gamma$ acts simply transitively
on the vertices of $P_{96}$
and this polytope therefore admits a Cayley structure.
The set $H$ of generators has 9 elements;
we explain how to obtain it.
We begin by identifying vertices of $P_{96}$
with oriented edges of $\{3,4,3\}$ as in the construction above.
Let $p_0$ be the vertex $(1,(1+i+j+k)/2)$:
by simple geometric considerations, its nine neighbours are
$(1,(1+i-j-k)/2)$, $(1,(1-i+j-k)/2)$, $(1,(1-i-j+k)/2)$,
$((1+i+j-k)/2,1)$, $((1+i-j+k)/2,1)$, $((1-i+j+k)/2,1)$,
$((1+i+j+k)/2, (1+i+j-k)/2)$, $((1+i+j+k)/2, (1+i-j+k)/2)$ and
$((1+i+j+k)/2, (1-i+j+k)/2)$.
The (unique) elements of $\Gamma$
taking $p_0$ to its nine neighbours are
$(i,i)$, $(j,j)$, $(k,k)$,
$((1+i-j+k)/2,k)$, $((1-i+j+k)/2,j)$, $((1+i+j-k)/2,i)$,
$((-1+i+j-k)/2,i)$, $((-1+i-j+k)/2,k)$ and $((-1-i+j+k)/2,j)$.
Thus, if we choose $p_0$ to be the identity
for the Cayley structure, the nine elements of $\Gamma$ above
are the nine elements of $H$.
From the previous section, we have explicit matrices
for the 19 irreducible representations of $\Gamma$;
their dimensions are at most 4.
The spectrum of the adjacency matrix $X$ can now
be computed as in section 4.

\bigbreak
{\nobigbf 8. $P_{720}$}

\nobreak\bigskip\nobreak

Following Coxeter ([C], sections 8.1 and 8.9),
we take for vertices of $P_{720} = \left\{{3 \atop 3,5}\right\}$
the midpoints of the edges of the regular polytope $\{3,3,5\}$.
Each vertex of $P_{720}$ is surrounded by two icosahedra
and five octahedra; its vertex figure is a pentagonal prism.

We now describe the group $G_{14400} \subset O(4)$
of all isometries of $P_{720}$.
The group $G_{7200}$ of orientation preserving isometries
of $\{3,3,5\}$ has order $120 \times 60$,
since it is transitive on the 120 vertices
and the subgroup of such isometries keeping a given vertex fixed
equals the group of orientation preserving isometries
of the vertex figure, an icosahedron.
Thus, $G_{7200} = (Q_{120} \times Q_{120})/(-1,-1)$.
The group $G_{14400}$ is generated by $G_{7200}$
together with a reflection on a hyperplane preserving
the vertices of $\{3,3,5\}$.

Unfortunately, the technique of the previous sections
does not apply directly.

{\nobf Proposition:}
{\sl The polytope $P_{720}$ admits no Cayley structure.}

{\nobf Proof:}
Let $G_{720}$ be an arbitrary subgroup of order 720 of $G_{14400}$:
we prove that $G_{720}$ does not
act simply transitively on edges of $\{3,3,5\}$.
Let $G_{25}$ be a 5-Sylow subgroup of $G_{720}$.
Clearly, $G_{25}$ is contained in
$G_{7200} = (Q_{120} \times Q_{120})/(-1,-1)$
and, by lifting,
we obtain a 5-Sylow subgroup of $Q_{120} \times Q_{120}$
which is conjugate to
$\{1,q,q^2,q^3,q^4\} \times \{1,q,q^2,q^3,q^4\}$,
where $q$ is some quaternion of order 5.
Thus, $G_{720}$ contains some element $g$ conjugate to $(q,q)$,
whose action keeps some vertex $v \in \{3,3,5\}$ fixed
(since $(q,q)$ does).
The element $g$ permutes the 12 neighbours of $v$,
splitting them into orbits of size 1 or 5;
hence, there is a neighbour $w$ of $v$, and hence an edge $vw$,
which are kept fixed under $g$, as we wanted to show.
\qed

Lov\'asz describes a method ([L]) to reduce the problem
of computing the spectrum of a graph with a transitive group
$\Gamma$ of isomorphisms to each irreducible decomposition of $\Gamma$,
which could be applied to this example.
As before, instead of counting paths,
we prefer to work with the matrices for representations
in a modified version of Lov\'asz's technique.
For simplicity, we start by applying the procedure
to $\langle 3,5,3,5\rangle$, which also admits no Cayley structure,
since its isometry group, $A_5 \oplus \ZZ/(2)$,
has no subgroup of order 30 (the number of vertices).

As illustrated in Figure 8.1, 
the vertices of $\langle 3,5,3,5\rangle$ are the midpoints
of edges of a dodecahedron.
If instead of taking midpoints of edges we take two
suitably spaced points per edge, we obtain $\langle 3,10,10\rangle$,
which admits a Cayley structure described in section 5.
Recall that $a = (12345)$, $b = (253)$ and $H = \{b,b^{-1},ab\}$.
Edges between decagons correspond to the generator $ab = (ab)^{-1}$.
For each vertex of $\langle 3,5,3,5\rangle$ there are
two elements of $A_5$ (counting the identity) which keep it fixed.

\midinsert\smallskip

\centerline{\psfig{file=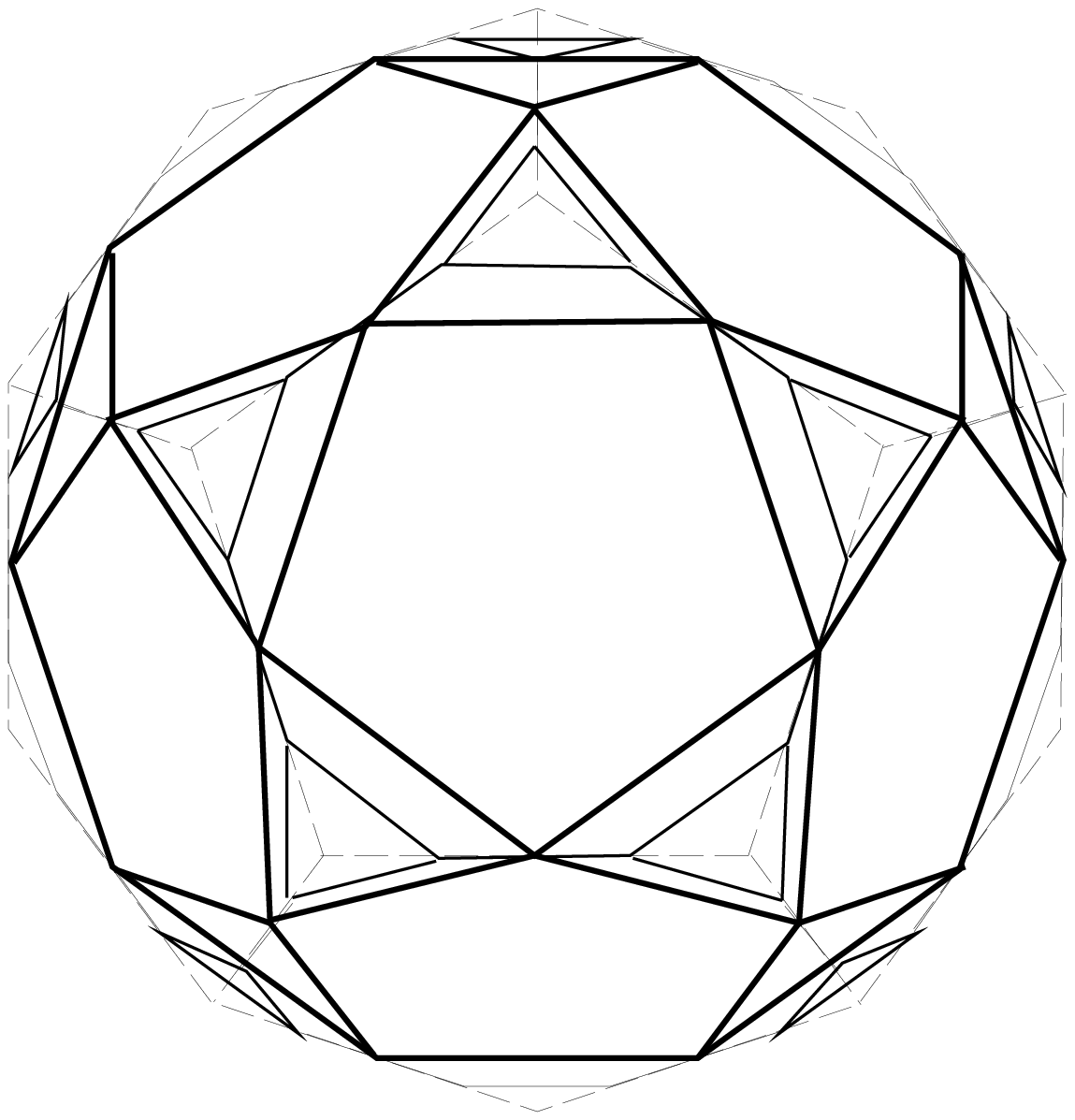,width=3.0in}}

\smallskip

\centerline{\eightbf Figure 8.1}

\smallskip\endinsert

Define a linear injection
$A_1: V_{\langle 3,5,3,5\rangle} \to V_{\langle 3,10,10\rangle}$
such that the value of $A_1(v)$ at a vertex $p$
of $\langle 3,10,10\rangle$
is the value of $v$ at the midpoint of the only edge
between decagons containing $p$
(recall that $V_P$ is the set of complex valued
functions on the vertices of the polytope $P$).
Conversely, define
$A_2: V_{\langle 3,10,10\rangle} \to V_{\langle 3,5,3,5\rangle}$
so that the value of $A_2(w)$ at a vertex $q$ of
$\langle 3,5,3,5\rangle$ is the sum of the values of $w$
at the two ends of the edge of $\langle 3,10,10\rangle$ containing $q$.
Thus, $A_2 A_1 = 2I$ and $A_1 A_2 = I + R_{ab}$,
where $R$ is the right multiplication action.
We claim that
$$X_{3535} = A_2(R_b + R_{b^{-1}})A_1;$$
Figure 8.2 illustrates the equality
at a basis vector of $V_{\langle 3,5,3,5\rangle}$.
The matrix $Y = (R_b + R_{b^{-1}})A_1A_2 =
(R_b + R_{b^{-1}})(I + R_{ab})$ has the same spectrum
as $X_{3535}$, up to 30 extra zero eigenvalues.
It is now clear that $Y$ splits into the irreducible representations
and its spectrum is computed in the usual manner.

\midinsert\smallskip

\centerline{\psfig{file=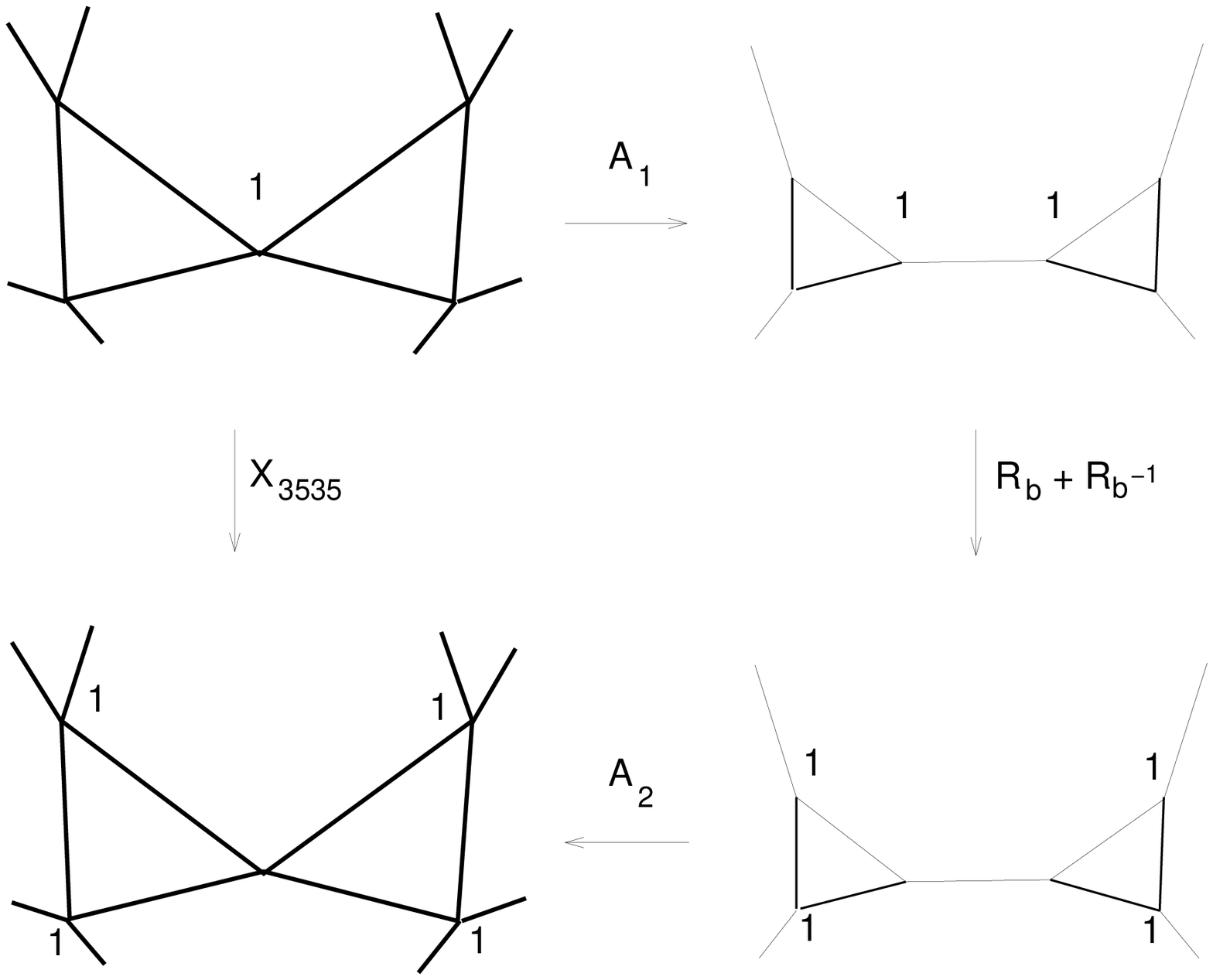,width=3.0in}}

\smallskip

\centerline{\eightbf Figure 8.2}

\smallskip\endinsert

We are ready to consider $P_{720}$.
As with $\langle 3,5,3,5\rangle$,
take two points in each edge of $\{3,3,5\}$
to obtain a (non-semi-regular) polytope $P_{1440}$ with 1440 vertices.
Call edges of $P_{1440}$
contained in edges of $\{3,3,5\}$ {\sl special}.
The group $G_{7200}$ of orientation preserving
isometries of $\{3,3,5\}$ does not act simply transitively
on the vertices of $P_{1440}$.
Happily, its subgroup $G_{1440} = (Q_{120} \times Q_{24})/(-1,-1)$ does;
in other words, $P_{1440}$ admits a Cayley structure.
Indeed, the first factor $Q_{120}$ guarantees that
$G_{1440}$ acts transitively on the vertices of $\{3,3,5\}$.
The subgroup of $G_{1440}$ keeping the vertex 1 fixed
consists of the 12 elements of the form $(q,q)$, where $q \in Q_{24}$.
These act simply transitively on the 12 neighbours in $\{3,3,5\}$
of the vertex 1,
as can be checked using the coordinate system in Section 6.
The identity for the Cayley structure of $P_{1440}$
is chosen to be the vertex
between 1 and $(\tau + i - \taub j)/2$ which is closer to 1.
Using coordinates and quaternion multiplication,
the reader may check that
$H = \{g_0,g_1,g_2,g_3,g_4,g_s\}$, where
$g_0 = (i,i)$, $g_1 = ((1+i+j+k)/2,(1+i+j+k)/2)$,
$g_2 = ((1+i+j-k)/2,(1+i+j-k)/2)$,
$g_3 = ((1-i-j+k)/2,(1-i-j+k)/2)$,
$g_4 = ((1-i-j-k)/2,(1-i-j-k)/2)$ and
$g_s = ((-\taub i - j + \tau k)/2, k)$;
special edges correspond to $g_s$.
Notice that $g_0^{-1} = g_0$,
$g_1^{-1} = g_4$, $g_2^{-1} = g_3$ and $g_s^{-1} = g_s$.

Again, define linear transformations
$A_1: V_{P_{720}} \to V_{P_{1440}}$ and
$A_2: V_{P_{1440}} \to V_{P_{720}}$:
$A_1(v)$ at a vertex $p$ of $P_{1440}$
is the value of $v$ at the midpoint of the special edge containing $p$
and $A_2(w)$ at a vertex $q$ of
$P_{720}$ is the sum of the values of $w$
at the two ends of the special edge containing $q$.
Thus, $A_2 A_1 = 2I$ and $A_1 A_2 = I + R_{g_s}$.
Also, $$X_{P_{720}} =
A_2(R_{g_0} + R_{g_1} + R_{g_2} + R_{g_3} + R_{g_4})A_1.$$
In order to prove this equality,
we relate the adjacencies of $P_{720}$ and $P_{1440}$.
A vertex $p$ of $P_{720}$ has 10 neighbours
and is the midpoint of a special edge of $P_{1440}$
with vertices $q$ and $q'$.
The vertex $q$ has six neighbours: $qg_0,\ldots,qg_4$ and $q' = qg_s$.
Similarly, the six neighbours of $q'$ are $q'g_0,\ldots,q'g_4$ and
$q = q'g_s$.
Omitting the repetitions of $q$ and $q'$,
the special edges containing the remaining 10 points
have as midpoints the 10 neighbours $p_0,\ldots,p_9$ of $p$.
The process $p \mapsto \{q,q'\} \mapsto
\{qg_0,\ldots,qg_4,q'g_0,\ldots,q'g_4\} \mapsto
\{p_0,\ldots,p_9\}$ corresponds the successive application on 
a basis vector of $V_{720}$ of the transformations
$A_1$, $R_{g_0} + \cdots + R_{g_4}$ and $A_2$,
finishing the proof of the equality.
The rest is routine by now: 
$Y = (R_{g_0} + \cdots + R_{g_4})(I + R_{g_s})$
has the same spectrum as $X_{P_{720}}$,
up to 720 extra zero eigenvalues.
Finally, split $Y$ into the 32 irreducible representations
of $G_{1440}$ to compute its spectrum.

\bigbreak
\noindent{\bigbf References}

\nobreak\bigskip\nobreak

\parindent=50pt

\item{[B]}Babai, L.,
{\sl Spectra of Cayley graphs},
J. Comb. Theory, Ser. B, 27 (1979), 180-189.

\item{[BB]}Blind, G. and Blind, R.,
{\sl The semiregular polytopes},
Comment. Math. Helvetici, 66 (1991), 150-154.

\item{[BCN]}Brouwer, A. E., Cohen, A. M. and Neumaier, A.,
{\sl Distance-regular graphs},
Springer-Verlag, New York, 1989.

\item{[C]}Coxeter, H. S. M.,
{\sl Regular polytopes},
Dover, New York, 1973.

\item{[C1]}Coxeter, H. S. M.,
{\sl Regular and semi-regular polytopes, I},
Math. Zeitschrift, 46 (1940), 380-407.

\item{[C2]}Coxeter, H. S. M.,
{\sl Regular and semi-regular polytopes, II},
Math. Zeitschrift, 188 (1985), 559-591.

\item{[C3]}Coxeter, H. S. M.,
{\sl Regular and semi-regular polytopes, III},
Math. Zeitschrift, 200 (1988), 3-45.

\item{[CCNPW]}Conway, J. H., Curtis, R. T., Norton, S. P.,
Parker, R. A. and Wilson, R. A.,
{\sl Atlas of finite groups},
Clarendon Press, Oxford, 1985.

\item{[CDS]}Cvetkovi\'c, D. M., Doob, M. and Sachs, H.,
{\sl Spectra of graphs},
Academic Press, New York, 1979.

\item{[CS]}Chung, F. R. K. and Sternberg, S.,
{\sl Laplacian and vibrational spectra for homogeneous graphs},
J. Graph Theory, 16 (1992), 605-627.

\item{[Ga]}Gantmacher, F.,
{\sl The theory of matrices},
Chelsea, New York, 1974.

\item{[Go]}Gosset, T.,
{\sl On the regular and semi-regular figures
in space of $n$ dimensions},
Messenger of Math., 29 (1900), 43-48.

\item{[K]}Kepler, J.,
{\sl Harmonice Mundi},
Opera Omnia, vol. 5, Frankfurt, 1864.

\item{[L]}Lov\'asz, L., 
{\sl Spectra of graphs with transitive groups},
Periodica Math. Hungarica, 6 (2) (1975), 191-195.

\item{[MKS]}Magnus, W., Karrass, A. and Solitar, D.,
{\sl Combinatorial group theory},
Interscience, New York, 1966.

\item{[MT]}Mneimn\'e, R. and Testard, F.,
{\sl Introduction \`a la th\'eorie des groupes de Lie classiques},
Hermann, Paris, 1986.

\item{[PS]}Petersdorf, M. and Sachs, H.,
{\sl Spektrum und Automorphismengruppe eines Graphen},
Combinatorial Theory and its Applications
(Colloq. Math. Soc. J. Bolyai, 4),
Amsterdam--London, 1970, 891-907.

\item{[ST]}Saldanha, N. and Tomei, C.,
{\sl Spectra of regular polytopes},
Discrete Comput. Geom., 7 (1992), 403-414.

\item{[Sc]}Schl\"afli, L.,
{\sl Theorie der vielfachen Kontinuit\"at},
Denkschriften der Schweizerischen naturforschenden Gesellschaft,
38 (1901), 1-237.

\item{[Se]}Serre, J. P.,
{\sl Repr\'esentations lin\'eaires des groupes finis},
Hermann, Paris, 1978.

\item{[Su]}Sugiura, M.,
{\sl Unitary representations an harmonic analysis, an introduction},
second edition, North-Holland--Kodansha, Tokyo, 1990.

\bigskip
\bigbreak

\vbox{\obeylines
\parskip=0pt
\parindent=0pt

Nicolau C. Saldanha, PUC-Rio and IMPA
nicolau@impa.br; http://www.impa.br/$\sim$nicolau/

\smallskip

Carlos Tomei, PUC-Rio and IMPA
tomei@impa.br

\bigskip

Departamento de Matem\'atica, PUC-Rio
Rua Marqu\^es de S\~ao Vicente, 225
Rio de Janeiro, RJ 22453-900, Brasil

\smallskip

IMPA
Estr. Dona Castorina, 110
Rio de Janeiro, RJ 22460-320, Brasil
}

\bye